\documentclass[11pt]{article}
\usepackage{stmaryrd}
\usepackage{tipa}
\usepackage{pifont}
\usepackage{amssymb}
\usepackage{latexsym}
\usepackage{eucal}
\usepackage{amsthm,mathtools}
\usepackage{amsmath}
\usepackage{mathabx}
\usepackage[dvips]{graphics}
\usepackage{epsfig}
\newtheorem{theorem}{Theorem}[section]
\newtheorem{proposition}[theorem]{Proposition}
\newtheorem{lemma}[theorem]{Lemma}

\newtheorem{corollary}[theorem]{Corollary}

\begin{document}
\title{Connectivity of Double Graphs}%
\author{Zhayida Simayijiang and Elkin Vumar \thanks{Corresponding author.
Email, vumar@xju.edu.cn.}\\[1pt]
College of Mathematics and System Sciences\\[1pt]
Xinjiang University, Urumqi 830046, P.R.China\\[1pt]}
\date{}
\maketitle 
\begin{abstract}

The double graph of a graph $G$ is defined as $\mathcal{D}[G]$ = $G
\times T_2$, where \(T_2\) is the total graph with 2 vertices and
$\times$ stands for the Kronecker product of graphs. In this paper,
sufficient conditions for double graphs to be maximum
vertex-connected, maximum edge-connected are presented.
\\
\\
$Keywords$: Double graph;maximum vertex-connected;maximum
edge-connected
\end{abstract}
\bigskip
\section{Introduction}

\bigskip

In this paper, unless specified otherwise, we consider only finite
simple graphs (i.e., without loops and multiple edges). As usual
$V(G)$ and $E(G)$ denote the sets of vertices and edges of $G$,
respectively, and $adj$ denotes the adjacency relation of $G$. A
vertex of degree $1$ in a graph is called a leaf vertex(or simply, a
leaf), and an edge incident with a leaf is called a leaf edge. For
notation and terminology not defined here we refer to West
\cite{West}.

The total graph $T_n$ on $n$ vertices is the graph obtained from the
complete graph $K_n$ by adding a loop to every vertex. The double
graph of a graph $G$ is defined as $\mathcal{D}[G]$ = $G \times
T_2$, where \(T_2\) is the total graph with 2 vertices, and $\times$
stands for the Kronecker product of graphs. The Kronecker product
$G\times H$ of two graphs $G$ and $H$ is the graph  with $V(G\times
H)=V(G)\times V(H)$ and with adjacency defined by $(u_1,v_1)$ $adj$
$(u_2,v_2)$ if and only if $u_1$ $adj$ $u_2$ in $G$ and $v_1$ $adj$
$v_2$ in H. In \cite{Emanuele} it was observed that there is a kind
of general construction which can be performed on every simple
graph. The class of double graphs with this construction turned out
to have several interesting properties. Some known results on double
graphs are given in \cite{Emanuele}.

If $V(T_2)=\{0,1\}$, then $G_0=\{(v,0): v\in V(G)\}$ and
$G_1=\{(v,1): v\in V(G)\}$ induce two subgraphs of $\mathcal{D}[G]$
both isomorphic to $G$ such that ${G_0}\bigcap{G_1}=\emptyset$ and
${G_0}\bigcup{G_1}$ induces a spanning subgraph of $\mathcal{D}[G]$.
We call $\{G_0,G_1\}$ the canonical decomposition of
$\mathcal{D}[G]$.

As a generalization of double graphs, we define ${\mathcal{D}_n}[G]$
= $G \times T_n$, where \(T_n\) is the total graph with $n$
vertices. Similarly, we call $\{G_0,G_1,...,G_{n-1}\}$ the canonical
decomposition of $\mathcal{D}_n[G]$. Note that
${\mathcal{D}_2}[G]=\mathcal{D}[G]$.

In what follows, for a graph $G=(V,E)$, we use $p(G)$ and $q(G)$ (or
simply $p$ and $q$) to denote $|V|$ and $|E|$, respectively. A graph
$G=(V,E)$ is maximum vertex-connected (in short, max-$\kappa$) if
$\kappa(G)$=$[\frac{2q(G)}{p(G)}]$, where $\kappa(G)$ is the
vertex-connectivity of $G$. Similarly, $G$ is maximum edge-connected
(in short, max-$\lambda$) if $\lambda(G)$=$[\frac{2q(G)}{p(G)}]$,
where $\lambda(G)$ is the edge-connectivity of $G$. Different
sufficient conditions for a graph to be max-$\kappa$ or
max-$\lambda$ have been recently given in the literature, see
Refs.[5-9].

\bigskip
\bigskip

We decided to write this paper as some graphical parameters of
double graphs that perhaps deserve to be better known. In Section 2,
we consider max-$\kappa$ of double graphs and in Section 3, we
consider max-$\lambda$ of double graphs.

\begin{proposition}\cite{Emanuele}
$\kappa(\mathcal{D}[G])=2\kappa(G)$.\qed
\end{proposition}

\begin{lemma}
$p({\mathcal{D}_n}[G])=np(G)$, $q({\mathcal{D}_n}[G])=n^2q(G)$,
$deg_{\mathcal{D}_n[G]}(u,v)=ndeg_{G}(u)$.\qed
\end{lemma}

In [1], some basic properties of double graphs $\mathcal{D}[G]$ are
given, it is not difficult to extend some of them to
$\mathcal{D}_n[G]$, here we list them in below.
\begin{proposition}
For any graph $G\ne K_1$ the following properties hold.\\
1. $G$ is connected if and only if $\mathcal{D}_n[G]$ is connected.\\
2. If $G$ is connected, then every pair of vertices of
$\mathcal{D}_n[G]$ belongs to a cycle.\\
3. Every edge of $\mathcal{D}_n[G]$ belongs to a $2n$-cycle.\\
4. In $\mathcal{D}_n[G]$ there is neither cut vertex nor cut
edge.\\
5. If $G$ is connected, then $\mathcal{D}_n[G]$ is a block. \qed
\end{proposition}
\begin{proposition}
For any graph $G$, $\mathcal{D}_n[G]$ is bipartite if and only if
 $G$ is bipartite. \qed
\end{proposition}
\begin{proposition}
For any graph $G\ne K_1$ the following traversability properties
hold.\\
1. Let $G$ be a connected graph, then $\mathcal{D}_n[G]$ is
eulerian if and only if  $G$ is eulerian or $n$ is even.\\
2. If $G$ is Hamiltonian, then so is $\mathcal{D}_n[G]$.
\begin{proof}
The proof of 1 is obvious, so we give the proof of 2. Let
$\{G_0,G_1,...,G_{n-1} \}$ be the canonical decomposition of
$\mathcal{D}_n[G]$. Let $\gamma$ be a spanning cycle of $G$, $uv$
and $u'v'$ be edges of $\gamma$ which are not incident with, and
$\gamma'$ be the path obtained from $\gamma$ by removing the edge
$uv$, $\pi$ and $\eta$ are the two components of
$\gamma-\{uv,u'v'\}$. Let $\gamma_0'$ and $\gamma_{n-1}'$ be the
corresponding paths of $\gamma-\{uv\}$ in $G_0$ and $G_{n-1}$,
respectively. Moreover, let $\pi_i=v_i...v'_i$ and
$\eta_i=u_i...u'_i$ be respectively the corresponding paths of $\pi$
and $\eta$ in $G_i$, for $i=1,2,...,n-2$.
 Then  \begin{eqnarray*}
\begin{split}
C=&\gamma_0' \bigcup \{v_0u_1\}\bigcup \eta_1 \bigcup \{u'_1v'_2\}
\bigcup
\pi_2 \bigcup ...\bigcup\gamma_{n-1}' \bigcup%
\\
&\{v_{n-1}u_{n-2}\}\bigcup\eta_{n-2}...\bigcup\pi_1\bigcup
\{v_1u_0\}
\end{split}
\end{eqnarray*} is a spanning cycle of $\mathcal{D}_n[G]$.(See
Fig.1).

\end{proof}
\end{proposition}

 \begin{center}
\begin{picture}(80,80)
\put(-90,40){\circle*{4}}%
\put(-90,-10){\circle*{4}}%
\put(-60,40){\circle*{4}}%
\put(-60,-10){\circle*{4}}%
\put(-30,40){\circle*{4}}%
\put(-30,-10){\circle*{4}}%
\put(0,40){\circle*{4}}%
\put(0,-10){\circle*{4}}%
\put(30,40){\circle*{4}}%
\put(30,-10){\circle*{4}}%
\put(40,15){\circle*{2}}%
\put(45,15){\circle*{2}}%
\put(50,15){\circle*{2}}%
\put(55,15){\circle*{2}}%
\put(60,15){\circle*{2}}%
\put(65,15){\circle*{2}}%

\put(80,40){\circle*{4}}%
\put(80,-10){\circle*{4}}%
\put(110,40){\circle*{4}}%
\put(110,-10){\circle*{4}}%
\put(140,40){\circle*{4}}%
\put(140,-10){\circle*{4}}%

\put(-142,12){$\gamma_0'$}%
\put(182,12){$\gamma_{n-1}'$}%
\put(-45,85){$\eta_1$}%
\put(-45,-60){$\pi_1$}%
\put(15,85){$\eta_2$}%
\put(15,-60){$\pi_2$}%
\put(95,85){$\eta_{n-2}$}%
\put(95,-60){$\pi_{n-2}$}%
\put(-88,42){$u_0$}%
\put(-88,-15){$v_0$}%
\put(-58,42){$u_1$}%
\put(-58,-15){$v_1$}%
\put(-28,42){$u_1'$}%
\put(-28,-15){$v_1'$}%
\put(2,42){$u_2'$}%
\put(2,-15){$v_2'$}%
\put(30,42){$u_2$}%
\put(30,-15){$v_2$}%
\put(58,42){$u_{n-2}'$}%
\put(58,-15){$v_{n-2}'$}%
\put(112,42){$u_{n-2}$}%
\put(112,-15){$v_{n-2}$}%
\put(142,42){$u_{n-1}$}%
\put(142,-15){$v_{n-1}$}%
\qbezier(-90,40)(-105,120)(-120,40)%
\qbezier(-90,-10)(-105,-90)(-120,-10)%
\qbezier(-120,40)(-125,15)(-120,-10)%

\qbezier(-90,40)(-75,15)(-60,-10)%
\qbezier(-90,-10)(-75,15)(-60,40)

\qbezier(-60,40)(-45,120)(-30,40)%
\qbezier(-60,-10)(-45,-90)(-30,-10)%

\qbezier(-30,40)(-15,15)(0,-10)%
\qbezier(-30,-10)(-15,15)(0,40)

\qbezier(0,40)(15,120)(30,40)%
\qbezier(0,-10)(15,-90)(30,-10)%

\qbezier(80,40)(95,120)(110,40)%
\qbezier(80,-10)(95,-90)(110,-10)%

\qbezier(110,40)(125,15)(140,-10)%
\qbezier(110,-10)(125,15)(140,40)%

\qbezier(140,40)(155,120)(170,40)%
\qbezier(140,-10)(155,-90)(170,-10)%
\qbezier(170,40)(175,15)(170,-10)%
\put(-10,-80){$Fig. 1$}

\end{picture}
\end{center}\vskip 3 cm


\begin{proposition}
$\kappa({\mathcal{D}_n}[G])=n\kappa(G)$.
\begin{proof}
Let $S$ be a minimum vertex cut of ${\mathcal{D}_n}[G]$. The sets
$S_i=S\bigcap V(G_i)$, $i=0,1,...,n-1$ are vertex cuts of
$G_0,G_1,...,G_{n-1}$, respectively. Then $|S_i|\ge\kappa(G)$ and
hence
$\kappa({\mathcal{D}_n}[G])\ge n\kappa(G)$.\\
Conversely, let $S$ be a vertex cut of $G$ and $S_i$ be the
corresponding sets in $G_i$, respectively, $i=0,1,...,n-1$. Then
$S_0\bigcup S_1\bigcup...\bigcup S_{n-1}$ is a vertex cut of
${\mathcal{D}_n}[G]$ and hence $\kappa({\mathcal{D}_n}[G])\le
n\kappa(G)$.
\end{proof}
\end{proposition}
\section{A sufficient condition on max-$\kappa$ of double graphs }

\bigskip

As simple examples, we can see that both $C_k$ and
$\mathcal{D}[C_k]$ are max-$\kappa$, while $K_{1,k-1}$ and $P_k$ are
max-$\kappa$, but $\mathcal{D}[K_{1,k-1}]$ and $\mathcal{D}[P_k]$
are not.
%
%
%
%
\begin{proposition}
Let $G$ be max-$\kappa$ and $q=tp+{t_0}$ with $0\le{t_0}\le{p-1}$.
If $0\le{t_0}<\frac{p}{4}$ or $\frac{p}{2}\le{t_0}<\frac{3p}{4}$,
then $\mathcal{D}[G]$ is max-$\kappa$.
\begin{proof}
By Lemma 1.1 and the definition of max-$\kappa$,
 we have $\kappa(\mathcal{D}[G])=2\kappa(G)=2[\frac{2q}{p}]$. On the
other hand, by Lemma 1.2,
$[\frac{2q(\mathcal{D}[G])}{p(\mathcal{D}[G])}]=[2\frac{2q}{p}]$.
Note that $q=tp+{t_0}$, $0\le{t_0}\le{p-1}$, we have\\

$[\frac{2q}{p}]=[\frac{2tp+2t_0}{p}]=[2t+\frac{2t_0}{p}]=\left\{
\begin{array}{ll}
2t, &0\le{t_0}<\frac{p}{2}\\
2t+1, &\frac{p}{2}\le{t_0}\le{p-1}
\end{array} \right. $

$[\frac{2q(\mathcal{D}[G])}{p(\mathcal{D}[G])}]=[2(2t+\frac{2t_0}{p})]=[4t+\frac{4t_0}{p}]$.\\
When $0\le{t_0}<\frac{p}{4}$, we have $[\frac{2q}{p}]=2t=\kappa(G)$
and
$[\frac{2q(\mathcal{D}[G])}{p(\mathcal{D}[G])}]=4t=\kappa(\mathcal{D}[G])$,
i.e., $\mathcal{D}[G]$ is max-$\kappa$. When
$\frac{p}{2}\le{t_0}<\frac{3p}{4}$, we have
$[\frac{2q}{p}]=2t+1=\kappa(G)$ and
$[\frac{2q(\mathcal{D}[G])}{p(\mathcal{D}[G])}]=4t+2=\kappa(\mathcal{D}[G])$,
i.e., $\mathcal{D}[G]$ is max-$\kappa$. It is easy to see that for
other value of $t_0$, $\mathcal{D}[G]$ is not max-$\kappa$.
\end{proof}
\end{proposition}
One may ask that if $\mathcal{D}[G]$ is max-$\kappa$ when $G$ is
not. The answer to this is negative as shown in the following.
\begin{proposition}
If $G$ is not max-$\kappa$, then $\mathcal{D}[G]$ is not
max-$\kappa$.
\begin{proof}
Suppose that $G$ is not max-$\kappa$ and $\mathcal{D}[G]$ is
max-$\kappa$. Then $\kappa(G)\ne[\frac{2q}{p}]$ and
$2\kappa(G)=\kappa(\mathcal{D}[G])=[2\frac{2q}{p}]\ne2[\frac{2q}{p}]$.
From the above inequalities, it is not difficult to deduce that
$\frac{p}{4}\le{t_0}<\frac{p}{2}$ or $\frac{3p}{4}\le{t_0}\le{p-1}$.
But when $\frac{p}{4}\le{t_0}<\frac{p}{2}$, we have
$[\frac{2q}{p}]=2t$ and
$[\frac{2q(\mathcal{D}[G])}{p(\mathcal{D}[G])}]=4t+1$, contradicting
the fact that
$[\frac{2q(\mathcal{D}[G])}{p(\mathcal{D}[G])}]=2\kappa(G)$ is even.
Similarly, when $\frac{3p}{4}\le{t_0}\le{p-1}$, we have
$[\frac{2q}{p}]=2t+1$ and
$[\frac{2q(\mathcal{D}[G])}{p(\mathcal{D}[G])}]=4t+3$, again a
contradiction. Hence $\mathcal{D}[G]$ is not max-$\kappa$, if $G$ is
not.
\end{proof}
\end{proposition}
\begin{theorem}
$\mathcal{D}[G]$ is max-$\kappa$ if and only if $G$ is max-$\kappa$
with $0\le{t_0}<\frac{p}{4}$ or $\frac{p}{2}\le{t_0}<\frac{3p}{4}$,
where $q=tp+{t_0}$, $0\le{t_0}\le{p-1}$.\qed
\end{theorem}

By using a similar argument, the result on max-$\kappa$ of double
graphs can easily be extended to graphs $\mathcal{D}_n[G]$. Hence we
have the following theorem on max-$\kappa$ of $\mathcal{D}_n[G]$.

\begin{theorem}
$\mathcal{D}_n[G]$ is max-$\kappa$ if and only if $G$ is
max-$\kappa$ with $0\le{t_0}<\frac{p}{2n}$ or $\frac{p}{2}\le {t_0}<
\frac{p(n+1)}{2n}$, where $q=tp+{t_0}$, $0\le{t_0}\le{p-1}$.\qed
\end{theorem}


%
%
%
\section{A sufficient condition on max-$\lambda$ of double graphs}
We start this section with some simple observations.

{\noindent}{\bf Fact1.} If $G$ is connected, then $\mathcal{D}[G]$
has no cut
edge. Consequently, $\lambda(\mathcal{D}[G])\ge2$. \\
\vskip 0.03cm {\noindent}{\bf Fact2.} If a connected graph $G$ has a
leaf vertex, then $\lambda(\mathcal{D}[G])=2$. In particular,
$\lambda(\mathcal{D}[T])=2$ for a tree $T$.
\begin{proposition}
If $G$ is a connected graph, then $\lambda(\mathcal{D}[G])\ge
2\lambda(G)$.
\begin{proof}
For any edge set $W\subseteq E(\mathcal{D}[G])$ with
$|W|<2\lambda(G)$, we need to show that $\mathcal{D}[G]-W$ is
connected. Suppose $W$ is such a set of $E(\mathcal{D}[G])$,
 $W=W_0\bigcup W_1\bigcup W_2$, where $W_0$ and $W_1$ are the corresponding
 edges of $W$ in $E(G_0)$ and $E(G_1)$, respectively.
 Let $R=\mathcal{D}[G]/[E(G_0)\bigcup E(G_1)]\cong G\times K_2$,
 $E(R)=\{(v_i,0)(v_j,1)|v_iv_j\in E(G)\}$, $W_2=W\bigcap E(R)$.
Since $R$ contains all vertices of $\mathcal{D}[G]$, we deduce that
if $R-W_2$ is connected, then so is $\mathcal{D}[G]-W$. Without loss
of
generality, we may assume $R-W_2$ is disconnected.\\
 {\bf Case 1.} $|W_2|\ge \lambda(G)$.\\
 In this case, we have $|W_0|< \lambda(G)$ and $|W_1|< \lambda(G)$, i.e., $G_0-W_0$ and
 $G_1-W_1$ are
 connected. Since $|E(R)|=2|E(G)|$, we have $|W_2|< 2\lambda(G)\le 2|E(G)|=|E(R)|$, and then in
 $E(R)-W_2$
  there is at leat one edge connecting $G_0$ and $G_1$, so $\mathcal{D}[G]-W$ is
 connected.\\
 {\bf Case2.} $|W_2|< \lambda(G)$.\\
 If $|W_0|<\lambda(G)$ and $|W_1|<\lambda(G)$, then both $G_0-W_0$ and
 $G_1-W_1$ are connected. Since $|W_2|<\lambda(G)$, as in Case 1,
 there is at least one edge in $E(R)-W_2$ connecting $G_0-W_0$ and
 $G_1-W_1$, and consequently $G-W$ is connected.\\
 Now assume that $|W_0|\ge \lambda(G)$ or $|W_1|\ge \lambda(G)$, say
 the former, then $|W_1|<\lambda(G)$ and $G_1-W_1$ is connected. If
 $G_0-W_0$ is connected, then we are done, hence assume $G_0-W_0$ is
 disconnected.\\
 Suppose, to the contrary, $G-W$ is disconnected, and $G'_1$ is a
 component of $G-W$. Then $G'_1$ is a component of $G_0-W_0$, since
 $G_1-W_1$ is connected. Since $G_0\cong G$, in $G_0$ there are at
 least $\lambda(G)$ edges between $V(G'_1)$ and $V(G_0)-V(G'_1)$. By
 the definition of $\mathcal{D}[G]$, there are at least $\lambda(G)$
 edges between $V(G'_1)$ and $V(G_1)$, and therefore $|W_2|\ge
 \lambda(G)$, a contradiction. Hence $G-W$ is connected and the
 proof is complete.
\end{proof}
\end{proposition}
Since $\delta(\mathcal{D}[G])=2\delta(G)$, we have the following
corollary.

\begin{corollary}
If $\lambda(G)=\delta(G)$, therefore $\lambda(\mathcal{D}[G])=
2\lambda(G)$.\qed
\end{corollary}
\begin{proposition}
If $G$ is a connected graph, then $\lambda(\mathcal{D}[G])=\left\{
\begin{array}{ll}
4\lambda(G), & if \lambda(G)\le \frac{\delta(G)}{2}\\
2\delta(G), & if \frac{\delta(G)}{2}<\lambda(G)<\delta(G)
\end{array} \right.$ \qed
\begin{proof}
Let $S$ be a minimum edge  cut of $G$, $S_0$ and $S_1$ be the
corresponding copies of $S$ in $G_0$, $G_1$, respectively. Set
$S_2=\{(u,0)(v,1),(v,0)(u,1)|e=uv\in S \}$, thus $|S_2|=2|S|$. Then
$S_0\bigcup S_1\bigcup S_2$ is an edge cut of $\mathcal{D}[G]$ and
hence $\lambda(\mathcal{D}[G])\le 4\lambda(G)$.\\
Let $W$ be an edge cut of $\mathcal{D}[G]$. The sets $W_0=W\bigcap
E(G_0)$, $W_1=W\bigcap E(G_1)$ and $W_2=W\bigcap E(R)$ are edge cuts
of $G_0$, $G_1$ and $R$, respectively. Then $|W_0|$,
$|W_1|\ge \lambda(G)$.\\
To get minimum edge cut of double graph, we consider $\lambda(G)$
and $\delta(G)$. There is two possible way to choose minimum edge
cut of $\mathcal{D}[G]$, one is choose from $G_0$, $G_1$ and $R$,
another one is choose all edges of $R$ as edge cut. So
$\lambda(\mathcal{D}[G])=min\{2\delta(G), 4\lambda(G),
2|E(G)|\}=min\{2\delta(G), 4\lambda(G)\}$. If $\lambda(G)\le
\frac{\delta(G)}{2}$, we have $4\lambda(G)\le 2\delta(G)$ . So
$\lambda(\mathcal{D}[G])= 4\lambda(G)$. If
$\frac{\delta(G)}{2}<\lambda(G)<\delta(G)$, we have
$2\delta(G)<4\lambda(G)<4\delta(G)$. In this case,
$\lambda(\mathcal{D}[G])=2\delta(G)$
\end{proof}
\end{proposition}
The following is an example to show that $\lambda(\mathcal{D}[G])=
4\lambda(G)$ when $\lambda(G)= \frac{\delta(G)}{2}$. In Fig.3 shows
$\lambda(\mathcal{D}[G])= 4\lambda(G)$ when $\lambda(G)<
\frac{\delta(G)}{2}$. In Fig.4 shows $\lambda(\mathcal{D}[G])=
2\delta(G)$ when $\frac{\delta(G)}{2}<\lambda(G)<\delta(G)$.\vskip
0.2 cm
  \begin{center}
\begin{picture}(60,80)

\put(-170,20){\circle*{4}}%
\put(-172,12){$u_1$}%
\put(-130,20){\circle*{4}}%
\put(-132,12){$u_3$}%
\put(-80,20){\circle*{4}} %
\put(-82,12){$u_4$}%
\put(-40,20){\circle*{4}}%
\put(-42,12){$u_6$}%
\put(-60,52){\circle*{4}} %
\put(-62,57){$u_5$}%
\put(-150,52){\circle*{4}} %
\put(-152,57){$u_2$}%

\put(-170,20){\line(1,0){130} }%
\put(-170,20){\line(3,5){20} }%
\put(-150,52){\line(3,-5){20} }%
\put(-80,20){\line(3,5){20} }%
\put(-60,52){\line(3,-5){20} }%

\put(-130,-30){$\lambda(G)=1,\delta(G)=2$}

\put(20,45){\circle*{4}} %
\put(60,45){\circle*{4}}%
\put(110,45){\circle*{4}} %
\put(150,45){\circle*{4}} %
\put(40,77){\circle*{4}} %
\put(130,77){\circle*{4}} %

\put(20,45){\line(1,0){130} }%

\put(20,45){\line(3,5){20} }%
\put(40,77){\line(3,-5){20} }%
\put(110,45){\line(3,5){20} }%
\put(130,77){\line(3,-5){20} }%

\put(20,20){\circle*{4}} %
\put(60,20){\circle*{4}}%
\put(110,20){\circle*{4}} %
\put(150,20){\circle*{4}} %
\put(40,-5){\circle*{4}} %
\put(130,-5){\circle*{4}} %

\put(20,20){\line(1,0){130} }%

\qbezier(20,45)(40,32.5)(60,20) %
\qbezier(60,45)(40,32.5)(20,20) %
\qbezier(60,45)(85,32.5)(110,20) %
\qbezier(110,45)(85,32.5)(60,20) %
\qbezier(110,45)(130,32.5)(150,20) %
\qbezier(150,45)(130,32.5)(110,20) %

\qbezier(20,45)(30,20)(40,-5) %
\qbezier(60,45)(50,20)(40,-5) %
\qbezier(110,45)(120,20)(130,-5) %
\qbezier(150,45)(140,20)(130,-5) %

\qbezier(20,20)(30,7.5)(40,-5) %
\qbezier(60,20)(50,7.5)(40,-5) %
\qbezier(110,20)(120,7.5)(130,-5) %
\qbezier(150,20)(140,7.5)(130,-5) %

\qbezier(40,77)(30,48.5)(20,20)%
\qbezier(40,77)(50,48.5)(60,20)%
\qbezier(130,77)(120,48.5)(110,20)%
\qbezier(130,77)(140,48.5)(150,20)%

\put(55,-30){$\lambda(\mathcal{D}[G])=4$}%
\put(-10,-60){$Fig. 2$}

\end{picture}
\end{center}\vskip 2 cm

  \begin{center}
\begin{picture}(60,80)

\put(-70,20){\circle*{4}}%

\put(-30,20){\circle*{4}}%

\put(20,20){\circle*{4}} %

\put(60,20){\circle*{4}}%

\put(-50,52){\circle*{4}} %
\put(-50,-12){\circle*{4}} %

\put(40,52){\circle*{4}} %
\put(40,-12){\circle*{4}} %

\put(-70,20){\line(1,0){130} }%
\put(-70,20){\line(3,5){20} }%
\put(-50,52){\line(3,-5){20} }%
\put(20,20){\line(3,5){20} }%
\put(40,52){\line(3,-5){20} }%

\qbezier(-70,20)(-60,4)(-50,-12) %
\qbezier(-30,20)(-40,4)(-50,-12)%

\qbezier(20,20)(30,4)(40,-12) %
\qbezier(60,20)(50,4)(40,-12)%

\qbezier(-50,52)(-50,20)(-50,-12) %
\qbezier(40,52)(40,21)(40,-12)

\put(-50,-40){$\lambda(G)=1, \delta(G)=3$}%
\put(-110,20){$G$}%
\put(-10,-60){$Fig. 3$}
\end{picture}
\end{center}\vskip 2 cm

  \begin{center}
\begin{picture}(60,80)

\put(-50,52){\circle*{4}} %
\put(-50,-12){\circle*{4}} %
\put(-50,20){\circle*{4}} %

\put(-80,39){\circle*{4}} %
\put(-80,1){\circle*{4}} %

\put(60,52){\circle*{4}} %
\put(60,-12){\circle*{4}}%
\put(60,20){\circle*{4}} %

\put(90,39){\circle*{4}} %
\put(90,1){\circle*{4}} %

\qbezier(-50,52)(-50,20)(-50,-12) %
\qbezier(-80,39)(-80,20)(-80,1) %
\qbezier(-80,39)(-65,45.5)(-50,52) %
\qbezier(-80,39)(-65,29.5)(-50,20) %
\qbezier(-80,39)(-65,13.5)(-50,-12) %
\qbezier(-80,1)(-65,26.5)(-50,52) %
\qbezier(-80,1)(-65,10.5)(-50,20) %
\qbezier(-80,1)(-65,-5.5)(-50,-12) %

\qbezier(-50,52)(5,52)(60,52) %
\qbezier(-50,20)(5,20)(60,20) %
\qbezier(-50,-12)(5,-12)(60,-12) %

\qbezier(60,52)(60,20)(60,-12) %
\qbezier(90,39)(90,20)(90,1) %
\qbezier(90,39)(75,45.5)(60,52) %
\qbezier(90,39)(75,29.5)(60,20) %
\qbezier(90,39)(75,13.5)(60,-12) %
\qbezier(90,1)(75,26.5)(60,52) %
\qbezier(90,1)(75,10.5)(60,20) %
\qbezier(90,1)(75,-5.5)(60,-12) %

\put(-70,-40){$\lambda(G)=3, \delta(G)=4, \lambda(\mathcal{D}[G])=
8$}%
\put(-110,20){$G$}%
\put(-10,-60){$Fig. 4$}
\end{picture}
\end{center}\vskip 2 cm

Corollary 3.2 and Proposition 3.3 yield the following theorem.
\begin{theorem}
$\lambda(\mathcal{D}[G])=\left\{
\begin{array}{ll}
2\lambda(G), & if \lambda(G)=\delta(G)\\
4\lambda(G), & if \lambda(G)\le \frac{\delta(G)}{2}\\
2\delta(G), & if \frac{\delta(G)}{2}<\lambda(G)<\delta(G)
\end{array} \right.$ \qed
\end{theorem}
\begin{proposition}
Let $G$ be max-$\lambda$ and $q=tp+{t_0}$ with $0\le{t_0}\le{p-1}$.
If $0\le{t_0}<\frac{p}{4}$ or $\frac{p}{2}\le{t_0}<\frac{3p}{4}$,
then $\mathcal{D}[G]$ is max-$\lambda$ when
$\lambda(\mathcal{D}[G])=2\lambda(G)$.
\begin{proof}
By Theorem 3.4 and the definition of max-$\lambda$,
 we have $\lambda(\mathcal{D}[G])=2\lambda(G)=2[\frac{2q}{p}]$. On the
other hand, by Lemma 1.2,
$[\frac{2q(\mathcal{D}[G])}{p(\mathcal{D}[G])}]=[2\frac{2q}{p}]$.
Note that $q=tp+{t_0}$, $0\le{t_0}\le{p-1}$, we have\\

$[\frac{2q}{p}]=[\frac{2tp+2t_0}{p}]=[2t+\frac{2t_0}{p}]=\left\{
\begin{array}{ll}
2t, &0\le{t_0}<\frac{p}{2}\\
2t+1, &\frac{p}{2}\le{t_0}\le{p-1}
\end{array} \right. $

$[\frac{2q(\mathcal{D}[G])}{p(\mathcal{D}[G])}]=[2(2t+\frac{2t_0}{p})]=[4t+\frac{4t_0}{p}]$.\\
When $0\le{t_0}<\frac{p}{4}$, we have $[\frac{2q}{p}]=2t=\lambda(G)$
and
$[\frac{2q(\mathcal{D}[G])}{p(\mathcal{D}[G])}]=4t=\lambda(\mathcal{D}[G])$,
i.e., $\mathcal{D}[G]$ is max-$\lambda$. When
$\frac{p}{2}\le{t_0}<\frac{3p}{4}$, we have
$[\frac{2q}{p}]=2t+1=\lambda(G)$ and
$[\frac{2q(\mathcal{D}[G])}{p(\mathcal{D}[G])}]=4t+2=\lambda(\mathcal{D}[G])$,
i.e., $\mathcal{D}[G]$ is max-$\kappa$. It is easy to see that for
other value of $t_0$, $\mathcal{D}[G]$ is not max-$\lambda$.
\end{proof}
\end{proposition}
\begin{proposition}
Let $G$ be max-$\lambda$, then $\mathcal{D}[G]$ is not max-$\lambda$
when $\lambda(\mathcal{D}[G])=4\lambda(G)$.
\end{proposition}
A natural question is that if $\mathcal{D}[G]$ is max-$\lambda$ when
$G$ is not. The answer is also negative.
\begin{proposition}
If $G$ is not max-$\lambda$, then $\mathcal{D}[G]$ is not
max-$\lambda$.
\begin{proof}
Suppose that $G$ is not max-$\lambda$ and $\mathcal{D}[G]$ is
max-$\lambda$. Then $\lambda(G)\ne[\frac{2q}{p}]$ and
$2\kappa(G)=\kappa(\mathcal{D}[G])=[2\frac{2q}{p}]\ne2[\frac{2q}{p}]$.
From the above inequalities, it is not difficult to deduce that
$\frac{p}{4}\le{t_0}<\frac{p}{2}$ or $\frac{3p}{4}\le{t_0}\le{p-1}$.
But when $\frac{p}{4}\le{t_0}<\frac{p}{2}$, we have
$[\frac{2q}{p}]=2t$ and
$[\frac{2q(\mathcal{D}[G])}{p(\mathcal{D}[G])}]=4t+1$, contradicting
the fact that
$[\frac{2q(\mathcal{D}[G])}{p(\mathcal{D}[G])}]=2\lambda(G)$ is
even. Similarly, when $\frac{3p}{4}\le{t_0}\le{p-1}$, we have
$[\frac{2q}{p}]=2t+1$ and
$[\frac{2q(\mathcal{D}[G])}{p(\mathcal{D}[G])}]=4t+3$, again a
contradiction. Hence $\mathcal{D}[G]$ is not max-$\lambda$, if $G$
is not.
\end{proof}
\end{proposition}
\begin{theorem}
$\mathcal{D}[G]$ is max-$\lambda$ when
$\lambda(\mathcal{D}[G])=2\lambda(G)$ if and only if $G$ is
max-$\lambda$ with $0\le{t_0}<\frac{p}{4}$ or
$\frac{p}{2}\le{t_0}<\frac{3p}{4}$, where $q=tp+{t_0}$,
$0\le{t_0}\le{p-1}$.\qed
\end{theorem}

It is not clear if $\mathcal{D}[G]$ is max-$\lambda$ when
$\frac{\delta(G)}{2}<\lambda(G)<\delta(G)$.

\begin{theorem}
$\lambda(\mathcal{D}_n[G])=\left\{
\begin{array}{ll}
n\lambda(G), & if \lambda(G)=\delta(G)\\
n^2\lambda(G), & if \lambda(G)\le \frac{\delta(G)}{2}\\
2n\delta(G), & if \frac{\delta(G)}{2}<\lambda(G)<\delta(G)
\end{array} \right.$ \qed
\end{theorem}

\begin{theorem}
$\mathcal{D}_n[G]$ is max-$\lambda$ if and only if $G$ is
max-$\lambda$ with $0\le{t_0}<\frac{p}{2n}$ or $\frac{p}{2}\le
{t_0}< \frac{p(n+1)}{2n}$, where $q=tp+{t_0}$,
$0\le{t_0}\le{p-1}$.\qed
\end{theorem}

\end{document}